\documentstyle[amsfonts,amsmath,12pt]{article}
\newtheorem{th}{Theorem}[section]

\newtheorem{cor}[th]{Corollary}
\newtheorem{defn}[th]{Definition}
\newenvironment{defn-new}{\begin{defn} \em}{\end{defn}}
\newtheorem{rem}[th]{Remark}
\newenvironment{rem-new}{\begin{rem} \em}{\end{rem}}
\newtheorem{ex}[th]{Example}
\newenvironment{ex-new}{\begin{ex} \em}{\end{ex}}

\newenvironment{notation-new}{\begin{rem} \em}{\end{rem}}

\newenvironment{agr-new}{\begin{rem} \em}{\end{rem}}

\makeatletter \@addtoreset{equation}{section} \makeatother

\makeatletter \@addtoreset{figure}{section} \makeatother

\begin{document}

\begin{center}
{\bf {\Large Clairaut anti-invariant submersion from nearly Kaehler manifold}%
}\\[0pt]

{\bf Punam Gupta and Amit Kumar Rai}
\end{center}

\bigskip

\noindent {\bf Abstract.} In the present paper, we investigate geometric
properties of Clairaut anti-invariant submersions whose total space is a
nearly Kaehler manifold. We obtain condition for Clairaut anti-invariant
submersion to be a totally geodesic map and also study Clairaut
anti-invariant submersions with totally umbilical fibers. In the last, we
introduce illustrative example.

\noindent {\bf 2010 Mathematics Subject Classification.}\medskip 53C12,
53C15, 53C20, 53C55


\section{Introduction}

Riemannian submersion between two Riemannian manifolds was first introduced
by O'Neill \cite{Neill} and Gray \cite{Gray}. After that Watson \cite{Watson}
introduced almost Hermitian submersions. Later, the notion of anti-invariant
submersions and Lagrangian submersion from almost Hermitian manifolds onto
Riemannian manifolds were introduced by Sahin \cite{Sahina} and studied by Ta%
\c{s}tan \cite{Tastan,Tastana}, G\"{u}nd\"{u}zalp \cite{Gundu}, Beri et al.
\cite{Beri}, Ali and Fatima \cite{Ali-Fatima}, in which the fibers of
submersion are anti-invariant with respect to the almost complex structure
of total manifold. After that several new types of Riemannian submersions
were defined and studied such as semi-invariant submersion \cite{Oz,Sahins},
slant submersion \cite{Erken,Sahinsl}, generic submersion \cite%
{Ali,Akyol,Chen,Ronsse}, hemi-slant submersion \cite{Tastanh}, semi-slant
submersion \cite{Park}, pointwise slant submersion \cite{Aykurt,Chen12,Lee5}
and conformal semi-slant submersion \cite{Akyolc}. Also, these kinds of
submersions were considered in different kinds of structures such as nearly
Kaehler, Kaehler, almost product, para-contact, Sasakian, Kenmotsu,
cosymplectic and etc. In book \cite{Sahinbook}, we find the recent
developments in this field.

In the theory of surfaces, Clairaut's theorem states that for any geodesic $%
\alpha $ on a surface of revolution $S$ , the function $r\sin \theta $ is
constant along $\alpha $, where $r$ is the distant from a point on the
surface to the rotation axis and $\theta $ is the angle between $\alpha $
and the meridian through $\alpha $. Bishop \cite{Bishop} introduced the idea
of Riemannian submersions and gave a necessary and sufficient conditions for
a Riemannian submersion to be Clairaut. Allison \cite{Allison} considered
Clairaut semi-Riemannian submersions and showed that such submersions have
interesting applications in the static space-times.

In \cite{Tastan1}, Tastan and Gerdan gave new Clairaut conditions for
anti-invariant submersions whose total manifolds are Sasakian and Kenmotsu
and got many interesting results. In \cite{Tastan-Aydin}, Tastan and Aydin
studied Clairaut anti-invariant submersions whose total manifolds are
cosymplectic. G\"{u}nd\"{u}zalp \ \cite{Gund} introduced Clairaut
anti-invariant submersions from a paracosymplectic manifold and gave
characterization theorems. In \cite{Lee1}, Lee et al. studied Clairaut
anti-invariant submersions whose total manifolds are Kaehler.

The geometrically interesting class of almost Hermitian manifolds is of
nearly Kaehler manifolds, which is one of the sixteen classes of almost
Hermitian manifolds and also obtained by Gray and Hervella in their
remarkable paper \cite{Gray-Her}. The geometrical meaning of nearly Kaehler
condition is that the geodesics on the manifolds are holomorphically planar
curves. Gray in \cite{Gray} studied nearly Kaehler manifolds broadly and
gave example of a non-Kaehlerian nearly Kaehler manifold, which is $6$%
-dimensional sphere.

Motivated by this, we study Clairaut anti-invariant submersions from nearly
Kaehler manifolds onto Riemannian manifolds. We also obtain conditions for
Clairaut Riemannian submersion to be totally geodesic map. We investigate
conditions for the Clairaut anti-invariant submersions to be a totally
umbilical map. Also, we provide some examples.

\section{Preliminaries}

An almost complex structure on a smooth manifold $M$ is a smooth tensor
field $\varphi $ of type $(1,1)$ such that $\varphi ^{2}=-I$. A smooth
manifold equipped with such an almost complex structure is called an almost
complex manifold. An almost complex manifold $\left( M,\varphi \right) $
endowed with a chosen Riemannian metric $g$ satisfying
\begin{equation}
g(\varphi X,\varphi Y)=g(X,Y)  \label{eq-ka1}
\end{equation}%
for all $X,Y\in TM$, is called an almost Hermitian manifold.

An almost Hermitian manifold $M$ is called a nearly Kaehler manifold \cite%
{Gray} if
\begin{equation}
\left( \nabla _{X}\varphi \right) Y+\left( \nabla _{Y}\varphi \right) X=0
\label{eq-ka2}
\end{equation}%
for all $X,Y\in TM$. If $\left( \nabla _{X}\varphi \right) Y=0$ for all $%
X,Y\in TM$, then $M$ is known as Kaehler manifold. Every Kaehler manifold is
nearly Kaehler but converse need not be true.

\begin{defn-new}
{\rm {\cite{Neill,Neill1}}}\ Let $(M,g_{m})$ and $(N,g_{n})$ be Riemannian
manifolds, where $\dim (M)=m$, $\dim (N)=n$ and $m>n$. A Riemannian
submersion $\pi :M\rightarrow N$ is a map of $M$ onto $N$ satisfying the
following axioms:

\begin{description}
\item[(i)] $\pi$ has maximal rank.

\item[(ii)] The differential $\pi _{\ast }$ preserves the lengths of
horizontal vectors.
\end{description}
\end{defn-new}

For each $q\in N$, $\pi ^{-1}(q)$ is an $(m-n)$-dimensional submanifold of $%
M $. The submanifolds $\pi ^{-1}(q)$, $q\in N$, are called fibers. A vector
field on $M$ is called vertical if it is always tangent to fibers. A vector
field on $M$ is called horizontal if it is always orthogonal to fibers. A
vector field $X$ on $M$ is called basic if $X$ is horizontal and $\pi $%
-related to a vector field $X^{\prime }$ on $N$, that is, $\pi _{\ast
}X_{p}=X_{\pi _{\ast }(p)}^{\prime }$ for all $p\in M.$ We denote the
projection morphisms on the distributions $\ker \pi _{\ast }$ and $(\ker \pi
_{\ast })^{\perp }$ by ${\cal V}$ and ${\cal H}$, respectively. The sections
of ${\cal V}$ and ${\cal H}$ are called the vertical vector fields and
horizontal vector fields, respectively. So
\[
{\cal V}_{p}=T_{p}\left( \pi ^{-1}(q)\right) ,\qquad {\cal H}%
_{p}=T_{p}\left( \pi ^{-1}(q)\right) ^{\perp }.
\]

The second fundamental tensors of all fibers $\pi ^{-1}(q),\ q\in N$ gives
rise to tensor field $T$ and $A$ in $M$ defined by O'Neill \cite{Neill} for
arbitrary vector field $E$ and $F$, which is
\begin{equation}
T_{E}F={\cal H}\nabla _{{\cal V}E}^{M}{\cal V}F+{\cal V}\nabla _{{\cal V}%
E}^{M}{\cal H}F,  \label{EQ2.9}
\end{equation}%
\begin{equation}
A_{E}F={\cal H}\nabla _{{\cal H}E}^{M}{\cal V}F+{\cal V}\nabla _{{\cal H}%
E}^{M}{\cal H}F,  \label{EQ2.8}
\end{equation}%
where ${\cal V}$ and ${\cal H}$ are the vertical and horizontal projections.

On the other hand, from equations (\ref{EQ2.9}) and (\ref{EQ2.8}), we have%
\begin{equation}
\nabla _{V}W=T_{V}W+\widehat{\nabla }_{V}W,  \label{EQ2.10}
\end{equation}%
\begin{equation}
\nabla _{V}X={\cal H}\nabla _{V}X+T_{V}X,  \label{EQ2.11}
\end{equation}%
\begin{equation}
\nabla _{X}V=A_{X}V+{\cal V}\nabla _{X}V,  \label{EQ2.12}
\end{equation}%
\begin{equation}
\nabla _{X}Y={\cal H}\nabla _{X}Y+A_{X}Y,  \label{EQ2.13}
\end{equation}%
for all $V,W\in \Gamma (\ker \pi _{\ast })$ and $X,Y\in \Gamma (\ker \pi
_{\ast })^{\perp },$ where ${\cal V}\nabla _{V}W=\widehat{\nabla }_{V}W.$ If
$X$ is basic, then $A_{X}V={\cal H}\nabla _{V}X.$

It is easily seen that for $p\in M,$ $U\in {\cal V}_{p}$ and $X\in {\cal H}%
_{p}$ the linear operators
\[
{T}_{U},{A}_{X}:T_{p}M\rightarrow T_{p}M
\]%
are skew-symmetric, that is,
\begin{equation}
g({A}_{X}E,F)=-g(E,{A}_{X}F)\text{ and }g({T}_{U}E,F)=-g(E,{T}_{U}F),
\label{EQ2.14}
\end{equation}%
for all $E,F\in $ $T_{p}M.$ We also see that the restriction of ${T}$ to the
vertical distribution ${T}|_{\ker \pi _{\ast }\times \ker \pi _{\ast }}$ is
exactly the second fundamental form of the fibres of $\pi $. Since ${T}_{U}$
is skew-symmetric, therefore $\pi $ has totally geodesic fibres if and only
if ${T}\equiv 0$.

Let $\pi :(M,g_{m})\rightarrow (N,g_{n})$ be a smooth map between Riemannian
manifolds. Then the differential $\pi _{\ast }$ of $\pi $ can be observed a
section of the bundle $Hom(TM,\pi ^{-1}TN)\rightarrow M$, where $\pi ^{-1}TN$
is the bundle which has fibres $\left( \pi ^{-1}TN\right) _{x}=T_{f(x)}N$
has a connection $\nabla $ induced from the Riemannian connection $\nabla
^{M}$ and the pullback connection. Then the second fundamental form of $\pi $
is given by
\begin{equation}
(\nabla \pi _{\ast })(E,F)=\nabla _{E}^{N}\pi _{\ast }F-\pi _{\ast }(\nabla
_{E}^{M}F),\text{ \ for all }E,F\in \Gamma (TM),  \label{EQ2.15}
\end{equation}%
where $\nabla ^{N}$ is the pullback connection (\cite{Baird,Bejancu}). We
also know that $\pi $ is said to be totally geodesic map \cite{Baird} if $%
(\nabla \pi _{\ast })(E,F)=0,$ for all $E,F\in \Gamma (TM)$.

Let $\pi $ be an anti-invariantquarians Riemannian submersion from
nearly Kaehler manifold $(M,\varphi ,g_{m})$ onto Riemannian
manifolds $(N,g_{n})$. For any arbitrary tangent vector fields $U$
and $V$ on $M$, we set
\begin{equation}
(\nabla _{U}\varphi )V=P_{U}V+Q_{U}V  \label{eq-c1}
\end{equation}%
where $P_{U}V,Q_{U}V$ denote the horizontal and vertical part of $(\nabla
_{U}\varphi )V$, respectively. Clearly, if $M$ is a Kaehler manifold then $%
P=Q=0$.

If $M$ is a nearly Kaehler manifold then $P$ and $Q$ satisfy
\begin{equation}
P_{U}V=-P_{V}U,\qquad Q_{U}V=-Q_{V}U.  \label{eq-c2}
\end{equation}%
Consider
\[
\left( \ker \pi _{\ast }\right) ^{\bot }=\varphi \ker \pi _{\ast }\oplus \mu
,
\]%
where $\mu $ is the complementary distribution to $\varphi \ker \pi _{\ast }$
in $(ker\pi _{\ast })^{\bot }$ and $\varphi \mu \subset \mu $.

For $X\in \Gamma (ker\pi _{\ast })^{\perp }$, we have
\begin{equation}
\varphi X=\alpha X+\beta X,  \label{EQ3.4}
\end{equation}%
where $\alpha X\in \Gamma (\ker \pi _{\ast })$ and $\beta X\in \Gamma (\mu )$%
. If $\mu =0$, then an anti-invariant submersion is known as Lagrangian
submersion.

\begin{defn-new}
\cite{Lee}Let $(M,\varphi,g)$ be an almost Hermitian manifold and $N$ be a
Riemannian manifold with Riemannian metric $g_{n}$. Suppose that there
exists a Riemannian submersion $\pi :M\rightarrow N$, such that the vertical
distribution $\ker \pi _{\ast }$ is anti-invariant with respect to $\varphi $%
, i.e., $\varphi \ker \pi _{\ast }\subseteq \ker \pi _{\ast }^{\bot }$ .
Then, the Riemannian submersion $\pi $ is called an anti-invariant
Riemannian submersion. We will briefly call such submersions as
anti-invariant submersions.
\end{defn-new}

Let $S$ be a revolution surface in ${\Bbb R}^{3}$ with rotation axis $L$.
For any $p\in S$, we denote by $r(p)$ the distance from $p$ to $L$. Given a
geodesic $\alpha :J\subset {\Bbb R}\rightarrow S$ on $S$, let $\theta (t)$
be the angle between $\alpha (t)$ and the meridian curve through $\alpha
(t),t\in I$. A well-known Clairaut's theorem says that for any geodesic on $%
S $, the product $r\sin \theta $ is constant along $\alpha $, i.e., it is
independent of $t$. In the theory of Riemannian submersions, Bishop \cite%
{Bishop} introduces the notion of Clairaut submersion in the following way.

\begin{defn-new}
\cite{Bishop} A Riemannian submersion $\pi :(M,g)\rightarrow (N,g_{n})$ is
called a Clairaut submersion if there exists a positive function $r$ on $M $%
, such that, for any geodesic $\alpha $ on $M$, the function $(r\circ \alpha
)\sin \theta $ is constant, where, for any $t,\ \theta (t)$ is the angle
between $\dot{\alpha}(t)$ and the horizontal space at $\alpha (t)$.
\end{defn-new}

He also gave the following necessary and sufficient condition for a
Riemannian submersion to be a Clairaut submersion:

\begin{th}
\label{th-bis} \cite{Bishop} Let $\pi :(M,g)\rightarrow (N,g_{n})$ be a
Riemannian submersion with connected fibers. Then, $\pi $ is a Clairaut
submersion with $r=e^{f}$ if and only if each fiber is totally umbilical and
has the mean curvature vector field $H=-{grad}f$, where ${grad}f$ is the
gradient of the function $f$ with respect to $g$.
\end{th}

\section{Anti-invariant Clairaut Submersions from nearly \newline
Kaehler Manifolds}

In this section, we give new Clairaut conditions for anti-invariant
submersions from nearly Kaehler manifolds after giving some auxiliary
results.

\begin{th}
\label{th1} Let $\pi $ be an anti-invariant submersion from a nearly Kaehler
manifold $(M,\varphi ,g)$ onto a Riemannian manifold $(N,g_{n})$. If $%
h:J\subset
\mathbb{R}
\rightarrow M$ is a regular curve and $U(s)$ and $X(s)$ are the vertical and
horizontal parts of the tangent vector field $\dot{h}(s)=W$ of $h(s)$,
respectively, then $h$ is a geodesic if and only if along $h$
\begin{equation}
A_{X}\varphi U+A_{X}\beta X+T_{U}\beta X+{\cal V}\nabla _{X}\alpha
X+T_{U}\varphi U+\hat{\nabla}_{U}\alpha X=0,  \label{eq-1.}
\end{equation}%
\begin{equation}
{\cal H}\left( \nabla _{\dot{h}}\varphi U+\nabla _{\dot{h}}\beta X\right)
+A_{X}\alpha X+T_{U}\alpha X=0.  \label{eq-2.}
\end{equation}
\end{th}

\noindent {\bf Proof.} Let $\pi $ be an anti-invariant submersion from a
nearly Kaehler manifold $(M,\varphi ,g)$ onto a Riemannian manifold $%
(N,g_{n})$. Since $\varphi ^{2}\dot{h}=-\dot{h}$. Taking the covariant
derivative of this and using (\ref{eq-ka2}), we have
\begin{equation}
\left( \nabla _{\dot{h}}\varphi \right) \varphi \dot{h}+\varphi \left(
\nabla _{\dot{h}}\varphi \dot{h}\right) =-\nabla _{\dot{h}}\dot{h}.
\label{eq-1}
\end{equation}%
Since $U(s)$ and $X(s)$ are the vertical and horizontal parts of the tangent
vector field $\dot{h}(s)=W$ of $h(s)$, that is, $\dot{h}=U+X$. So (\ref{eq-1}%
) becomes
\begin{eqnarray}
-\nabla _{\dot{h}}\dot{h} &=&\varphi \left( \nabla _{U+X}\varphi
(U+X)\right) +P_{\dot{h}}\varphi \dot{h}+Q_{\dot{h}}\varphi \dot{h}
\nonumber \\
&=&\varphi \left( \nabla _{U}\varphi U+\nabla _{X}\varphi U+\nabla
_{U}\varphi X+\nabla _{X}\varphi X\right) +P_{\dot{h}}\varphi \dot{h}+Q_{%
\dot{h}}\varphi \dot{h}  \nonumber \\
&=&\varphi \left( \nabla _{U}\varphi U+\nabla _{X}\varphi U+\nabla
_{U}\left( \alpha X+\beta X\right) +\nabla _{X}\left( \alpha X+\beta
X\right) \right)   \nonumber \\
&&+P_{\dot{h}}\varphi \dot{h}+Q_{\dot{h}}\varphi \dot{h}.  \label{eq-2}
\end{eqnarray}%
Using (\ref{EQ2.10})-(\ref{EQ2.13}) in (\ref{eq-2}), we get
\begin{eqnarray}
-\nabla _{\dot{h}}\dot{h} &=&\varphi \left( {\cal H}\left( \nabla _{\dot{h}%
}\varphi U+\nabla _{\dot{h}}\beta X\right) +A_{X}\alpha X+A_{X}\beta
X+A_{X}\varphi U\right.   \nonumber \\
&&\left. +T_{U}\beta X+T_{U}\alpha X+{\cal V}\nabla _{X}\alpha
X+T_{U}\varphi U+\hat{\nabla}_{U}\alpha X\right) +P_{\dot{h}}\varphi \dot{h}%
+Q_{\dot{h}}\varphi \dot{h}.  \label{eq-3}
\end{eqnarray}%
Since $\varphi ^{2}X=-X$, on differentiation, we have
\[
\varphi \left( \nabla _{Y}\varphi X\right) +\left( \nabla _{Y}\varphi
\right) \varphi X=-\nabla _{X}Y,
\]%
\[
\varphi ^{2}\left( \nabla _{Y}X\right) +\varphi \left( \nabla _{Y}\varphi
\right) X+\left( \nabla _{Y}\varphi \right) \varphi X=-\nabla _{X}Y,
\]%
using (\ref{eq-c1}) in above, we obtain
\begin{equation}
\varphi \left( P_{Y}X+Q_{Y}X\right) =-P_{Y}\varphi X-Q_{Y}\varphi X.
\label{eq-c3}
\end{equation}%
By (\ref{eq-c3}), we have
\[
\varphi \left( P_{\dot{h}}\varphi \dot{h}+Q_{\dot{h}}\varphi \dot{h}\right)
=P_{\dot{h}}\dot{h}+Q_{\dot{h}}\dot{h},
\]%
since $P$ and $Q$ are antisymmetric, so
\begin{equation}
\varphi \left( P_{\dot{h}}\varphi \dot{h}+Q_{\dot{h}}\varphi \dot{h}\right)
=0.  \label{eq-c4}
\end{equation}%
Using (\ref{eq-c4}) and equating the vertical and horizontal part of (\ref%
{eq-3}), we obtain
\[
{\cal V}\varphi \nabla _{\dot{h}}\dot{h}=A_{X}\varphi U+A_{X}\beta
X+T_{U}\beta X+{\cal V}\nabla _{X}\alpha X+T_{U}\varphi U+\hat{\nabla}%
_{U}\alpha X,
\]%
\[
{\cal H}\varphi \nabla _{\dot{h}}\dot{h}={\cal H}\left( \nabla _{\dot{h}%
}\varphi U+\nabla _{\dot{h}}\beta X\right) +A_{X}\alpha X+T_{U}\alpha X.
\]%
By using above equations we can say that $h$ is geodesic if and only if (\ref%
{eq-1.}) and (\ref{eq-2.}) hold.

\begin{th}
Let $\pi $ be an anti-invariant submersion from a nearly Kaehler manifold $%
(M,\varphi ,g)$ onto a Riemannian manifold $(N,g_{n})$. Also, let $%
h:J\subset {\Bbb R}\rightarrow M$ be a regular curve and $U(s)$ and $X(s)$
are the vertical and horizontal parts of the tangent vector field $\dot{h}%
(s)=W$ of $h(s)$. Then $\pi $\ is a Clairaut submersion with $r=e^{f}$ if
and only if along $h$
\[
g({grad}f,X)g(U,U)=g({\cal H}\nabla _{\dot{h}}\beta X+A_{X}\alpha
X+T_{U}\alpha X+P_{\dot{h}(s)}U,\varphi U).
\]
\end{th}

\noindent {\bf Proof.} Let $h:J\subset
\mathbb{R}
\rightarrow M$ be a geodesic on $M$ and $\ell =\left\Vert \dot{h}%
(s)\right\Vert ^{2}$. Let $\theta (s)$ be the angle between $\dot{h}(s)$ and
the horizontal space at $h(s)$. Then
\begin{equation}
g(X(s),X(s))=\ell \cos ^{2}\theta (s),  \label{eq-4}
\end{equation}%
\begin{equation}
g(U(s),U(s))=\ell \sin ^{2}\theta (s).  \label{eq-5}
\end{equation}%
Differentiating (\ref{eq-5}), we get
\begin{equation}
2g(\nabla _{\dot{h}(s)}U(s),U(s))=2\ell \sin \theta (s)\cos \theta (s)\frac{%
d\theta (s)}{ds}.  \label{eq-6.}
\end{equation}%
Using (\ref{eq-ka1}) in (\ref{eq-6.}), we have
\[
g({\cal H}\nabla _{\dot{h}(s)}\varphi U(s),\varphi U(s))-g((\nabla _{\dot{h}%
(s)}\varphi )U(s),\varphi U(s))=\ell \sin \theta (s)\cos \theta (s)\frac{%
d\theta (s)}{ds}.
\]%
Now by use of (\ref{eq-c1}), we have
\[
g({\cal H}\nabla _{\dot{h}(s)}\varphi U(s),\varphi U(s))-g(P_{\dot{h}%
(s)}U+Q_{\dot{h}(s)}U,\varphi U(s))=\ell \sin \theta (s)\cos \theta (s)\frac{%
d\theta (s)}{ds}
\]%
Along the curve $h$, using Theorem \ref{th1}, we obtain
\[
-g({\cal H}\nabla _{\dot{h}}\beta X+A_{X}\alpha X+T_{U}\alpha X+P_{\dot{h}%
(s)}U,\varphi U(s))=\ell \sin \theta (s)\cos \theta (s)\frac{d\theta (s)}{ds}%
.
\]%
Now, $\pi $ is a Clairaut submersion with $r=e^{f}$ if and only if $\frac{d}{%
ds}\left( e^{f}\sin \theta \right) =0$. Therefore
\begin{eqnarray*}
e^{f}\left( \frac{df}{ds}\sin \theta +\cos \theta \frac{d\theta }{ds}\right)
&=&0, \\
e^{f}\left( \frac{df}{ds}\ell \sin ^{2}\theta +\ell \sin \theta \cos \theta
\frac{d\theta }{ds}\right)  &=&0.
\end{eqnarray*}%
So, we obtain
\begin{equation}
\frac{df}{ds}(h(s))g(U(s),U(s))=g({\cal H}\nabla _{\dot{h}}\beta
X+A_{X}\alpha X+T_{U}\alpha X+P_{\dot{h}(s)}U,\varphi U(s)),  \label{eq-6}
\end{equation}%
Since $\frac{df}{ds}(h(s))=g({grad}f,\dot{h}(s))=g({grad}f,X)$. Therefore by
using (\ref{eq-6}), we get result.

\begin{th}
\label{th2} Let $\pi $ be an Clairaut anti-invariant submersion from a
nearly Kaehler manifold $(M,\varphi ,g)$ onto a Riemannian manifold $%
(N,g_{n})$ with $r=e^{f}$. Then
\[
A_{\varphi W}\varphi X+Q_{W}\varphi X=X(f)W
\]%
for $X\in \left( \ker \pi _{\ast }\right) ^{\bot }$, $W\in \ker \pi _{\ast }$
and $\varphi W$ is basic.
\end{th}

\noindent {\bf Proof.} Let $\pi $ be an anti-invariant submersion from a
nearly Kaehler manifold $(M,\varphi ,g)$ onto a Riemannian manifold $%
(N,g_{n})$ with $r=e^{f}$. We know that any fiber of Riemannian submersion $%
\pi $ is totally umbilical if and only if
\begin{equation}
T_{V}W=g(V,W)H,  \label{eq-7.}
\end{equation}%
for all $V,W\in \Gamma (\ker \pi _{\ast })$, where $H$ denotes the mean
curvature vector field of any fiber in $M$. By using Theorem \ref{th-bis}
and (\ref{eq-7}), we have
\begin{equation}
T_{V}W=-g(V,W){grad}f.  \label{eq-8.}
\end{equation}%
Let $X\in \mu $ and $V,W\in \Gamma (\ker \pi _{\ast })$, then by using (\ref%
{eq-ka1}) and (\ref{eq-ka2}), we have
\begin{equation}
g(\nabla _{V}\varphi W,\varphi X)=g(\varphi \nabla _{V}W+(\nabla _{V}\varphi
)W,\varphi X)=g(\nabla _{V}W,X)+g(P_{V}W+Q_{V}W,\varphi X).  \label{eq-c7}
\end{equation}%
By using (\ref{eq-ka1}), we have
\[
g(\varphi Y,Z)=-g(Y,\varphi Z),
\]%
taking covariant derivative of above, we get
\[
g(\left( \nabla _{X}\varphi \right) Y,Z)=-g(Y,\left( \nabla _{X}\varphi
\right) Z),
\]%
using (\ref{eq-c1}), we get
\begin{eqnarray}
g(P_{X}Y+Q_{X}Y,Z) &=&-g(Y,P_{X}Z+Q_{X}Z)  \nonumber \\
&=&g(Y,P_{Z}X+Q_{Z}X).  \label{eq-c5}
\end{eqnarray}%
Using (\ref{eq-c5}), we have
\begin{equation}
g(P_{W}\varphi X+Q_{W}\varphi X,V)=g(\varphi X,P_{V}W+Q_{V}W)  \label{eq-c6}
\end{equation}%
Using (\ref{EQ2.10}), (\ref{eq-8.}), (\ref{eq-c6}) in (\ref{eq-c7}), we have
\[
g(\nabla _{V}\varphi W,\varphi X)=-g(V,W)\left( {grad}f,X\right)
+g(V,Q_{W}\varphi X).
\]%
Since $\varphi W$ is basic, so ${\cal H}\nabla _{V}\varphi W=A_{\varphi W}V$%
, therefore we have
\[
g(A_{\varphi W}V,\varphi X)=-g(V,W)\left( {grad}f,X\right) +g(V,Q_{W}\varphi
X),
\]%
\begin{equation}
g(V,A_{\varphi W}\varphi X)+g(V,Q_{W}\varphi X)=g(V,W)\left( {grad}%
f,X\right)   \label{eq-13}
\end{equation}%
because $A$ is anti-symmetric. By using (\ref{eq-13}), we get result.

\begin{th}
\label{th3} Let $\pi $ be a Clairaut anti-invariant submersion from a nearly
Kaehler manifold $(M,\varphi ,g)$ onto a Riemannian manifold $(N,g_{n})$
with $r=e^{f}$ and ${grad}f\in \varphi \ker \pi _{\ast }$. Then either $f$
is constant on $\varphi \ker \pi _{\ast }$ or the fibres of $\pi $ are $1$%
-dimensional.
\end{th}

\noindent {\bf Proof.} Using (\ref{EQ2.10}) and (\ref{eq-8.}), we have
\[
g(\nabla _{V}W,\varphi U)=-g(V,W)g({grad}f,\varphi U),
\]%
where $U,V,W\in \Gamma (\ker \pi _{\ast })$. Since $g(W,\varphi U)=0$.
therefore we have
\begin{equation}
g(W,\nabla _{V}\varphi U)=g(V,W)g({grad}f,\varphi U).  \label{eq-9}
\end{equation}%
By use of (\ref{eq-ka1}) and (\ref{eq-c1}) in (\ref{eq-9}), we get
\[
g(W,Q_{V}U)-g(\varphi W,\nabla _{V}U)=g(V,W)g({grad}f,\varphi U).
\]%
By using (\ref{EQ2.10}), we obtain
\[
g(W,Q_{V}U)-g(\varphi W,T_{V}U)=g(V,W)g({grad}f,\varphi U).
\]%
Now, using (\ref{eq-8.}), we get
\begin{equation}
g(W,Q_{V}U)+g(V,U)g({grad}f,\varphi W)=g(V,W)g({grad}f,\varphi U)
\label{eq-10}
\end{equation}%
Take $V=U$ in (\ref{eq-10}), we have
\begin{equation}
g(V,V)g({grad}f,\varphi W)=g(V,W)g({grad}f,\varphi V).  \label{eq-11}
\end{equation}%
Take $V=U$ and interchange $V$ with $W$ in (\ref{eq-10}), we have
\begin{equation}
g(W,W)g({grad}f,\varphi V)=g(V,W)g({grad}f,\varphi W).  \label{eq-12}
\end{equation}%
By (\ref{eq-11}) and (\ref{eq-12}), we have
\[
g^{2}(V,W)g({grad}f,\varphi V)=g(V,V)g(W,W)g({grad}f,\varphi V).
\]%
Therefore either $f$ is constant on $\varphi \ker \pi _{\ast }$ or $V=aW$,
where $a$ is constant (by using Schwarz's Inequality for equality case).

\begin{cor}
Let $\pi $ be a Clairaut anti-invariant submersion from a nearly Kaehler
manifold $(M,\varphi ,g)$ onto a Riemannian manifold $(N,g_{n})$ with $%
r=e^{f}$ and ${grad}f\in \varphi \ker \pi _{\ast }$. If $\dim (\ker \pi
_{\ast })>1$, then the fibres of $\pi $ are totally geodesic if and only if $%
A_{\varphi W}\varphi X+Q_{W}\varphi X=0$ for $W\in \ker \pi _{\ast }$ such
that $\varphi W$ is basic and $X\in \mu $.
\end{cor}

\noindent {\bf Proof.} By Theorem \ref{th2} and Theorem \ref{th3}, we get
the result.

\begin{cor}
Let $\pi $ be an Clairaut Lagrangian submersion from a nearly Kaehler
manifold $(M,\varphi ,g)$ onto a Riemannian manifold $(N,g_{n})$ with $%
r=e^{f}$. Then either the fibres of $\pi $ are $1$-dimensional or they are
totally geodesic.
\end{cor}

\noindent {\bf Proof.} Let $\pi $ be an Clairaut Lagrangian submersion from
a Kaehler manifold $(M,\varphi ,g)$ onto a Riemannian manifold $(N,g_{n})$
with $r=e^{f}$, Then $\mu =\{0\}$. So $A_{\varphi W}\varphi X+Q_{W}\varphi
X=0$ always. \medskip

Lastly, we give some examples for Clairaut anti-invariant submersions from a
nearly Kaehler manifold. \medskip

\noindent {\bf Example.} Let $({\Bbb R}^{4},\varphi ,g)$ be a nearly Kaehler
manifold endowed with Euclidean metric $g$ on ${\Bbb R}^{4}$ given by

\[
g=\sum\limits_{i=1}^{4}dx_{i}^{2}
\]%
and canonical complex structure
\[
\varphi (x_{j})=\left\{
\begin{array}{cc}
-x_{j+1} & j=1,3 \\
x_{j-1} & j=2,4%
\end{array}%
\right. .
\]%
The $\varphi $-basis is $\left\{ e_{i}=\frac{\partial }{\partial x_{i}}%
|i=1,2,3,4\right\} $. Let $( {\Bbb R} ^{3},g_{1})$ be a Riemannian manifold
endowed with metric $g=\sum\limits_{i=1}^{3}dy_{i}^{2}$.

\noindent {\bf (i)} Consider a map $\pi :({\Bbb R}^{4},\varphi
,g)\rightarrow ({\Bbb R} ^{3},g_{1})$ defined by
\[
\pi (x_{1},x_{2},x_{3},x_{4})=\left( \frac{x_{1}+x_{2}}{\sqrt{2}}%
,x_{3},x_{4}\right) .
\]%
Then by direct calculations, we have
\[
\ker \pi _{\ast }=span\left\{ X_{1}=\left( \frac{\partial }{\partial x_{1}}-%
\frac{\partial }{\partial x_{2}}\right) \right\} ,
\]%
\[
\left( \ker \pi _{\ast }\right) ^{\bot }=span\left\{ X_{2}=\left( \frac{%
\partial }{\partial x_{1}}+\frac{\partial }{\partial x_{2}}\right) ,X_{3}=%
\frac{\partial }{\partial x_{3}},X_{4}=\frac{\partial }{\partial x_{4}}%
\right\}
\]%
and $\varphi X_{1}=-X_{2}$, therefore $\varphi \left( \ker \pi _{\ast
}\right) \subset \left( \ker \pi _{\ast }\right) ^{\bot }$. Thus, we can say
that $\pi $ is an anti-invariant Riemannian submersion. Since the fibers of $%
\pi $ are $1$-dimensional, therefore fibers are totally umbilical.

Consider the Koszul formula for Levi-Civita connection $\nabla $ for ${\Bbb R%
}^{4}$

\[
2g(\nabla _{X}Y,Z)=Xg(Y,Z)+Yg(Z,X)-Zg(X,Y)-g([Y,Z],X)-g([X,Z],Y)+g([X,Y],Z)
\]%
for all $X,Y,Z\in {\Bbb R}^{4}$. By simple calculations, we obtain
\[
\nabla _{e_{i}}e_{j}=0\quad \text{for all }i,j=1,2,3,4.
\]%
Hence $T_{X}Y=T_{Y}X=T_{X}X=0$ for all $X,Y\in \Gamma (\ker \pi _{\ast })$.
Therefore fibres of $\pi $ are totally geodesic. Thus $\pi $ is Clairaut
trivially.

\medskip

\noindent {\bf (ii)} Consider a map $\pi :({\ {\Bbb R}} ^{4},\varphi
,g)\rightarrow ({\ {\Bbb R}}^{3},g_{1})$ defined by
\[
\pi (x_{1},x_{2},x_{3},x_{4})=\left( \sqrt{x_{1}^{2}+x_{2}^{2}}%
,x_{3},x_{4}\right) .
\]%
Then by direct calculations, we have
\[
\ker \pi _{\ast }=span\left\{ X_{1}=\left( \frac{x_{2}}{\sqrt{%
x_{1}^{2}+x_{2}^{2}}}\frac{\partial }{\partial x_{1}}-\frac{x_{1}}{\sqrt{%
x_{1}^{2}+x_{2}^{2}}}\frac{\partial }{\partial x_{2}}\right) \right\} ,
\]%
\[
\left( \ker \pi _{\ast }\right) ^{\bot }=span\left\{ X_{2}=\left( \frac{x_{2}%
}{\sqrt{x_{1}^{2}+x_{2}^{2}}}\frac{\partial }{\partial x_{1}}+\frac{x_{1}}{%
\sqrt{x_{1}^{2}+x_{2}^{2}}}\frac{\partial }{\partial x_{2}}\right) ,X_{3}=%
\frac{\partial }{\partial x_{3}},X_{4}=\frac{\partial }{\partial x_{4}}%
\right\}
\]%
and $\varphi X_{1}=-X_{2}$, therefore $\varphi \left( \ker \pi _{\ast
}\right) \subset \left( \ker \pi _{\ast }\right) ^{\bot }$. Thus, we can say
that $\pi $ is an anti-invariant Riemannian submersion. Since the fibres of $%
\pi $ are $1$-dimensional, therefore fibres are totally umbilical. By using
Koszul formula, we obtain
\[
\nabla _{e_{i}}e_{j}=0\quad \text{for all }i,j=1,2,3,4.
\]%
Hence
\[
T_{X_{1}}X_{1}=-\left( \frac{x_{2}}{\sqrt{x_{1}^{2}+x_{2}^{2}}}\frac{%
\partial }{\partial x_{1}}+\frac{x_{1}}{\sqrt{x_{1}^{2}+x_{2}^{2}}}\frac{%
\partial }{\partial x_{2}}\right) .
\]%
Now, for the function $f=\ln (\sqrt{x_{1}^{2}+x_{2}^{2}})$ on $( {\Bbb R}%
^{4},\varphi ,g)$, the gradient of $f$ with respect to $g$ is given by
\[
{grad}f=\sum\limits_{i,j=1}^{4}g^{ij}\frac{\partial f}{\partial x_{i}}\frac{%
\partial }{\partial x_{j}}=\frac{x_{2}}{\sqrt{x_{1}^{2}+x_{2}^{2}}}\frac{%
\partial }{\partial x_{1}}+\frac{x_{1}}{\sqrt{x_{1}^{2}+x_{2}^{2}}}\frac{%
\partial }{\partial x_{2}}.
\]%
Therefore for $X_{1}\in \Gamma (\ker \pi _{\ast })$, $T_{X_{1}}X_{1}=-{grad}%
f $. Since $\left\Vert X_{1}\right\Vert =1$, so $T_{X_{1}}X_{1}=-\left\Vert
X_{1}\right\Vert ^{2}{grad}f$. By using theorem~\ref{th-bis}, we can say
that $\pi $ is an proper Clairaut anti-invariant submersion with $r=e^{f}$
for $f=\ln (\sqrt{x_{1}^{2}+x_{2}^{2}})$.

\noindent Department of Mathematics and Statistics \newline
Dr. Harisingh Gour University\newline
Sagar-470 003,Madhya Pradesh\newline
India\newline
Email:{\em pgupta\makeatletter @\makeatother dhsgsu.edu.in} \bigskip

\noindent University School of Basic and Applied Sciences\newline
Guru Gobind Singh Indraprastha University\newline Delhi, India
\\
Email:{\em rai.amit08au\makeatletter @\makeatother gmail.com}

\end{document}